\documentclass{article}

\title{\LARGE \textbf{On some Versions of Conjectures of Bondy and Jung }}
\author{Zh.G. Nikoghosyan\footnote{G.G. Nicoghossian (up to 1997)}}  

\begin{document}

\maketitle

\begin{abstract}

Using algebraic transformations and equivalent reformulations we derive a number of new results from some earlier ones (by the author) in more accepted terms closely related to well-known conjectures of Bondy and Jung including a number of classical results in hamiltonian graph theory (due to Dirac, Ore, Nash-Williams, Bondy, Jung and so on) as special cases.  A number of extended and strengthened versions of these conjectures are proposed. \\

\noindent\textbf{Key words}. Hamilton cycle, dominating cycle, long cycles, large cycles, Bondy's conjecture, Jung's Conjecture.  

\end{abstract}

\section{Introduction}

Most of classical results in hamiltonian graph theory are associated in conjectures of Bondy \cite{[3]} and Jung \cite{[9]} in terms of degree sums (minimum degree), connectivity conditions and special extreme cycles, called large cycles, with appropriate outside structures, including well-known Hamilton and dominating cycles as special cases with extremely simple outside structures.

Throughout this article we consider only finite undirected graphs without loops or multiple edges. A good reference for any undefined terms is \cite{[4]}. 

 The set of edges of a graph $G$ is denoted by $E(G)$. Paths and cycles in a graph $G$ are considered as subgraphs of $G$. If $Q$ is a path or a cycle, then the length of $Q$, denoted by $|Q|$, is $|E(Q)|$. Each vertex and edge in a graph can be interpreted as simple cycles of lengths 1 and 2, respectively. A graph $G$ is hamiltonian if $G$ contains a Hamilton cycle, i.e. a cycle containing every vertex of $G$.  For $C$ a longest cycle in $G$, let $\overline{p}$ and $\overline{c}$ denote the lengths of a longest path and a longest cycle in $G\backslash C$, respectively. We can suppose that $\overline{p}=-1$ when $C$ is a Hamilton cycle. 

Generally, a cycle $C$ in a graph $G$ is a large cycle if it dominates some certain subgraph structures  in $G$ in a sense that every such structure has a vertex in common with $C$. When $C$ dominates all vertices in $G$ (that is $\overline{p}=-1$) then $C$ is a Hamilton cycle. When $C$ dominates all edges in $G$ (that is $\overline{p}=0$) then $C$ is called a dominating cycle introduced by Nash-Williams \cite{[11]}. Further, if $C$ dominates all paths in $G$ of length at least some fixed integer $\lambda$ then $C$ is a $PD_\lambda$ (path dominating)-cycle. Finally, if $C$ dominates all cycles in $G$ of length at least  $\lambda$  then $C$ is a  $CD_\lambda$ (cycle dominating)-cycle. In particular, $PD_0$-cycles and $CD_1$-cycles are well-known Hamilton cycles and $PD_1$-cycles and $CD_2$-cycles are often called dominating cycles.  

Let $G$ be a graph of order $n$ and minimum degree $\delta$. The degree sum of $t$ smallest degrees among $t$ pairwise nonadjacent vertices will be denoted by $\sigma_t$.

In 1980, Bondy \cite{[3]} conjectured a common generalization of well-known theorems of  Ore \cite{[14]} (1960, $\lambda=1$) and Bondy \cite{[3]} (1980, $\lambda=2$).\\ 

\noindent\textbf{Conjecture A} \cite{[3]}. Let $G$ be a $\lambda$-connected $(\lambda\ge1)$ graph and $C$ a longest cycle in $G$. If 
$$
\frac{1}{\lambda+1}\sigma_{\lambda+1}\ge\frac{n+2}{\lambda+1}+\lambda-2
$$
then $\overline{p}\le\lambda-2$.\\

When $\lambda=3$, Conjecture A was proved in 1987 by Zou \cite{[15]}.

The minimum degree  version of Conjecture A contains two fundamental theorems on this subject due to Dirac \cite{[5]} (1952, $\lambda=1$) and Nash-Williams \cite{[11]} (1971, $\lambda=2$) as special cases.\\

\noindent\textbf{Conjecture B} \cite{[3]}.  Let  $G$ be a $\lambda$-connected $(\lambda\ge 1)$ graph and $C$ a longest cycle in $G$. If 
$$
\delta \ge\frac{n+2}{\lambda+1}+\lambda-2
$$
then $\overline{p}\le\lambda-2$.\\

For $\lambda=3$ Conjecture B was proved in 1981 by Jung \cite{[8]}.  

The first result related to Conjecture B was established in 2009 in terms of $CD_{\lambda}$-cycles.\\

\noindent\textbf{Theorem A} \cite{[12]}.  Let  $G$ be a $\lambda$-connected $(\lambda\ge 1)$ graph and $C$ a longest cycle in $G$. If 
$$
\delta \ge\frac{n+2}{\lambda+1}+\lambda-2
$$
then $C$ is a $CD_{\min\{\lambda,\delta-\lambda+1\}}$-cycle.\\

Theorem A can be reformulated in terms of Conjecture B conforming that the minimum degree $\overline{c}$-version of Bondy's conjecture is true with some strengthening.\\

\noindent\textbf{Theorem 1}.  Let  $G$ be a $\lambda$-connected $(\lambda\ge 1)$ graph and $C$ a longest cycle in $G$. If 
$$
\delta \ge\frac{n+2}{\lambda+1}+\lambda-2
$$
then $\overline{c}\le \min\{\lambda-1,\delta-\lambda\}$.\\

In this paper we prove two analogous strengthenings related to Conjecture B without any connectivity conditions.\\

\noindent\textbf{Theorem 2}.  Let  $G$ be a graph and $C$ a longest cycle in $G$. If 
$$
\delta \ge\frac{n+2}{\lambda+1}+\lambda-2
$$
then either $\overline{p}\le\min\{\lambda-2,\delta-\lambda-1\}$ or 
$\overline{p}\ge\max\{\lambda,\delta-\lambda+1\}$. \\

\noindent\textbf{Theorem 3}.  Let  $G$ be a graph and $C$ a longest cycle in $G$. If 
$$
\delta \ge\frac{n+2}{\lambda+1}+\lambda-2
$$
then either 
$\overline{c}\le\min\{\lambda-1,\delta-\lambda\}$ or $\overline{c}\ge\max\{\lambda+1,\delta-\lambda+2\}$. \\

In view of Theorems 1-3 Conjectures A and B can be naturally extended by adding $\overline{c}$-versions and can be essentially strengthened.\\ 

\noindent\textbf{Conjecture 1}.  Let $G$ be a $\lambda$-connected $(\lambda\ge1)$ graph and $C$ a longest cycle in $G$. If 
$$
\frac{1}{\lambda+1}\sigma_{\lambda+1}\ge\frac{n+2}{\lambda+1}+\lambda-2
$$
then 
$$
\overline{p}\le\min\Big\{\lambda-2,\frac{1}{\lambda+1}\sigma_{\lambda+1}-\lambda-1\Big\},\ \ \ \overline{c}\le\min\Big\{\lambda-1,\frac{1}{\lambda+1}\sigma_{\lambda+1}-\lambda\Big\}.
$$
\\

\noindent\textbf{Conjecture 2}.  Let  $G$ be a $\lambda$-connected $(\lambda\ge 1)$ graph and $C$ a longest cycle in $G$. If 
$$
\delta \ge\frac{n+2}{\lambda+1}+\lambda-2
$$
then $\overline{p}\le\min\{\lambda-2,\delta-\lambda-1\}$.\\

Now we turn to long cycle versions of above developments.

In 2001, Jung \cite{[9]} conjectured a common generalization of two fundamental theorems in hamiltonian graph theory due to Dirac \cite{[5]} (1952, $\lambda=2$) and  Jung \cite{[7]} (1978, $\lambda=3$).\\

\noindent\textbf{Conjecture C} \cite{[9]}.  Let  $G$ be a $\lambda$-connected $(\lambda\ge 1)$ graph and $C$ a longest cycle in $G$. If  $\overline{p}\ge\lambda-2$ then $|C|\ge \lambda(\delta-\lambda+2)$.\\

The degree sum version of Conjecture C containing the theorems of Bondy \cite{[2]} (1971, $\lambda=2$), Bermond \cite{[1]} (1976, $\lambda=2$), Linial \cite{[10]} (1976, $\lambda=2$), Fraisse and Jung \cite{[6]} (1989, $\lambda=3$) as special cases can be formulated as follows.\\

\noindent\textbf{Conjecture 3}.  Let $G$ be a  $\lambda$-connected $(\lambda\ge 2)$ graph and $C$ a longest cycle in $G$. If $\overline{p}\ge\lambda-2$ then 
$$
|C|\ge \lambda \Bigl(\frac{1}{\lambda}\sigma_{\lambda}-\lambda+2\Bigl).
$$
\\

The first result related to Conjecture C was established in 2009 in terms of $CD_{\lambda}$-cycles.\\

\noindent\textbf{Theorem B} \cite{[12]}.  Let $G$ be a  $(\lambda+1)$-connected $(\lambda\ge 0)$ graph and $C$ a longest cycle in $G$. Then either $|C|\ge (\lambda+1)(\delta-\lambda+1)$ or $C$ is a $CD_{\min\{\lambda,\delta-\lambda\}}$-cycle.\\

Theorem B can be reformulated in terms of Conjecture C conforming that the minimum degree $\overline{c}$-version of Jung's conjecture is true with some strengthening.\\

\noindent\textbf{Theorem 4}.  Let $G$ be a  $\lambda$-connected graph and $C$ a longest cycle in $G$. If $\overline{c}\ge\min\{\lambda-1,\delta-\lambda+1\}$ then $|C|\ge \lambda(\delta-\lambda+2)$.\\

In this paper we prove two analogous strengthenings related to Conjecture C without any connectivity conditions.   \\

\noindent\textbf{Theorem 5}.  Let  $G$ be a graph and $C$ a longest cycle in $G$. If  
$$
\min\{\lambda-2,\delta-\lambda\}\le \overline{p}\le\max\{\lambda-2,\delta-\lambda\}
$$
then $|C|\ge\lambda(\delta-\lambda+2)$.\\

\noindent\textbf{Theorem 6}.  Let  $G$ be a graph and $C$ a longest cycle in $G$. If  
$$
\min\{\lambda-1,\delta-\lambda+1\}\le \overline{c}\le\max\{\lambda-1,\delta-\lambda+1\}
$$
then $|C|\ge\lambda(\delta-\lambda+2)$.\\

In view of Theorem 4-6, conjectures C and 3 can be naturally extended by adding $\overline{c}$-versions and can be essentially strengthened.\\

\noindent\textbf{Conjecture 4}.  Let $G$ be a  $\lambda$-connected $(\lambda\ge1)$ graph and $C$ a longest cycle in $G$. If either 
$$
\overline{p}\ge\min\Big\{\lambda-2,\frac{1}{\lambda}\sigma_{\lambda}-\lambda\Big\} \ \mbox{or} \ \
\overline{c}\ge\min\Big\{\lambda-1,\frac{1}{\lambda}\sigma_{\lambda}-\lambda+1\Big\}
$$
then 
$$
|C|\ge \lambda \Bigl(\frac{1}{\lambda}\sigma_{\lambda}-\lambda+2\Bigr).
$$
\\

\noindent\textbf{Conjecture 5}.  Let  $G$ be a $\lambda$-connected $(\lambda\ge 1)$ graph and $C$ a longest cycle in $G$. If   $\overline{p}\ge\min\{\lambda-2,\delta-\lambda\}$  then $|C|\ge \lambda(\delta-\lambda+2)$.\\

To prove Theorems 2,3,5,6 we need the following two theorems.\\

\noindent\textbf{Theorem C} \cite{[13]} (1998).  Let $G$ be a graph and $C$ a longest cycle in $G$. Then $|C|\ge(\overline{p}+2)(\delta-\overline{p})$.  \\

\noindent\textbf{Theorem D} \cite{[13]} (2000).  Let $G$ be a graph and $C$ a longest cycle in $G$. Then $|C|\ge(\overline{c}+1)(\delta-\overline{c}+1)$.\\

\section{Proofs}

\noindent\textbf{Proofs of Theorems 1 and 4}. Theorem 1 follows immediately since it is not hard to see from the definition that $C$ is a $CD_{\lambda}$-cycle if and only if $\overline{c}\le \lambda-1$. To prove Theorem 4, we can reformulate Theorem B in the following way by taking $\lambda-1$ instead of $\lambda$: if $G$ is a $\lambda$-connected graph $(\lambda\ge 1)$ then either $|C|\ge \lambda(\delta-\lambda+2)$ or $C$ is a $CD_{\min\{\lambda-1,\delta-\lambda+1\}}$-cycle. Further, since $\overline{c}\ge \min\{\lambda-1,\delta-\lambda+1\}$ (by the hypothesis), we conclude that $C$ is not a $CD_{\min\{\lambda-1,\delta-\lambda+1\}}$-cycle, implying that $|C|\ge \lambda(\delta-\lambda+2)$.     \ \qquad     \rule{7pt}{6pt}      \\

\noindent\textbf{Proof of Theorem 2}. By the hypothesis, $n\le (\lambda+1)(\delta-\lambda+2)-2$. On the other hand, we have $n\ge|C|+\overline{p}+1$. Since $|C|\ge(\overline{p}+2)(\delta-\overline{p})$ (by Theorem C), we have 
$$
n\ge(\overline{p}+2)(\delta-\overline{p})+\overline{p}+1=(\overline{p}+2)(\delta-\overline{p}+1)-1.
$$
Thus
$$
(\lambda+1)(\delta-\lambda+2)\ge(\overline{p}+2)(\delta-\overline{p}+1)+1,
$$
which is equivalent to 
$$
(\lambda-\overline{p}-1)(\delta-\lambda-\overline{p})\ge1.
$$
Then we have either
$$
\lambda-\overline{p}-1\ge1 \ \  \mbox{and} \ \  \delta-\lambda-\overline{p}\ge1,
$$
which is equivalent to  $\overline{p}\le\min\{\lambda-2,\delta-\lambda-1\}$, or
$$
\lambda-\overline{p}-1\le -1 \ \  \mbox{and} \ \  \delta-\lambda-\overline{p}\le -1,
$$
which is equivalent to  $\overline{p}\ge\max\{\lambda,\delta-\lambda+1\}$.   \ \qquad     \rule{7pt}{6pt}\\

\noindent\textbf{Proof of Theorem 3}. By the hypothesis, $n\le (\lambda+1)(\delta-\lambda+2)-2$. On the other hand, we have $n\ge|C|+\overline{c}$. Since $|C|\ge(\overline{c}+1)(\delta-\overline{c}+1)$ (by Theorem D), we have 
$$
n\ge(\overline{c}+1)(\delta-\overline{c}+1)+\overline{c}=(\overline{c}+1)(\delta-\overline{c}+2)-1.
$$
Thus
$$
(\lambda+1)(\delta-\lambda+2)\ge(\overline{c}+1)(\delta-\overline{c}+2)+1,
$$
which is equivalent to 
$$
(\lambda-\overline{c})(\delta-\lambda-\overline{c}+1)\ge1.
$$
Then we have either
$$
\lambda-\overline{c}\ge1 \ \  \mbox{and} \ \  \delta-\lambda-\overline{c}+1\ge1,
$$
which is equivalent to  $\overline{c}\le\min\{\lambda-1,\delta-\lambda\}$, or we have
$$
\lambda-\overline{c}\le -1 \ \  \mbox{and} \ \  \delta-\lambda-\overline{c}+1\le -1,
$$
which is equivalent to  $\overline{c}\ge\max\{\lambda+1,\delta-\lambda+2\}$.        \ \qquad     \rule{7pt}{6pt}
\\

\noindent\textbf{Proof of Theorem 5}. We distinguish two cases.

\textbf{Case 1}. $\min\{\lambda-2,\delta-\lambda\}=\lambda-2$. 

By the hypothesis, $\lambda-2\le\overline{p}\le\delta-\lambda$. Then
$$
(\overline{p}-\lambda+2)(\delta-\overline{p}-\lambda)\ge0,
$$
which is equivalent to 
$$
(\overline{p}+2)(\delta-\overline{p})\ge\lambda(\delta-\lambda+2).
$$
Since $|C|\ge(\overline{p}+2)(\delta-\overline{p})$ (by Theorem C), we have $|C|\ge\lambda(\delta-\lambda+2)$.\\

\textbf{Case 2}. $\min\{\lambda-2,\delta-\lambda\}=\delta-\lambda$.

By the hypothesis, $\delta-\lambda\le \overline{p}\le \lambda-2$, implying that 
$$
(\overline{p}-\lambda+2)(\delta-\overline{p}-\lambda)\ge 0
$$
and we can argue as in Case 1.        \ \qquad     \rule{7pt}{6pt}   \\

\noindent\textbf{Proof of Theorem 6}. We distinguish two cases.

\textbf{Case 1}. $\min\{\lambda-1,\delta-\lambda+1\}=\lambda-1$. 

By the hypothesis, $\lambda-1\le\overline{c}\le\delta-\lambda+1$. Then
$$
(\overline{c}-\lambda+1)(\delta-\overline{c}-\lambda+1)\ge0,
$$
which is equivalent to 
$$
(\overline{c}+1)(\delta-\overline{c}+1)\ge\lambda(\delta-\lambda+2).
$$
Since $|C|\ge(\overline{c}+1)(\delta-\overline{c}+1)$ (by Theorem D), we have $|C|\ge\lambda(\delta-\lambda+2)$.\\

\textbf{Case 2}. $\min\{\lambda-1,\delta-\lambda+1\}=\delta-\lambda+1$.

By the hypothesis, $\delta-\lambda+1\le \overline{c}\le \lambda-1$, implying that 
$$
(\overline{c}-\lambda+1)(\delta-\overline{c}-\lambda+1)\ge 0
$$
and we can argue as in Case 1.         \ \qquad     \rule{7pt}{6pt}    \\

\noindent Institute for Informatics and Automation Problems\\ National Academy of Sciences\\
P. Sevak 1, Yerevan 0014, Armenia\\ 
E-mail: zhora@ipia.sci.am

\end{document}